\documentclass[11pt]{amsart}

\usepackage{color}
\usepackage{amsbsy}
\usepackage{amssymb,amsmath,amstext,mathrsfs}

\usepackage{color}
\usepackage{amsbsy}
\usepackage{amssymb,amsmath,amstext,mathrsfs}
\usepackage{graphicx}

\makeatletter
\def\verbatim{\interlinepenalty\@M \@verbatim
  \leftskip\@totalleftmargin\advance\leftskip2pc
  \frenchspacing\@vobeyspaces \@xverbatim}
\makeatother
\newtheorem{thm}{Theorem}[section]
\newtheorem{cor}[thm]{Corollary}
\newtheorem{lem}[thm]{Lemma}
\newtheorem{prop}[thm]{Proposition}

\theoremstyle{definition}

\theoremstyle{remark}
\newtheorem{rem}{Remark}[section]

\numberwithin{equation}{section}

\newtheorem{ex}[thm]{Example}

\newcommand{\pl}{\left(}
\newcommand{\pr}{\right)}

\newcommand{\ra}{\longrightarrow}

\newcommand{\cx}{\mathcal{A}}
\newcommand{\cy}{\mathcal{B}}

\newcommand{\dx}{C (X, I)}
\newcommand{\dy}{C (Y, I)}
\newcommand{\dz}{C (Z, I)}

\newcommand{\txy}{\varphi: \cx \ra \cy}

\newcommand{\fxy}{\varphi: \dx \ra \dy}

\newcommand{\nin}{\mathbb{N} \cup \{ \infty \} }
\newcommand{\ro}{\mathbb{R}_{\ge 0}}

\title[Multiplicative bijections]{Multiplicative bijections of semigroups of interval-valued continuous functions}

\author{Jes\'us Araujo}
\address{Departamento de Matem\'aticas,
Estad\'{\i}stica y Computaci\'on\\ Universidad de Cantabria\\
Facultad de Ciencias\\ Avda.
de los Castros, s. n.\\ E-39071 Santander, Spain}

\email{araujoj@unican.es}
\thanks{Research  partially supported by
the Spanish Ministry of Science and Education (Grant number MTM2006-14786).}

\date{}
\begin{document}

\begin{abstract}
We characterize all compact and Hausdorff spaces $X$ which satisfy that  for every  multiplicative bijection  $\varphi$ on $C(X, I)$, there exist a homeomorphism $\mu : X \ra X$ and a continuous map $p: X \ra (0, +\infty)$ such that
$$\varphi (f) (x) = f(\mu (x))^{p(x)}$$
for every $f \in \dx$ and $x \in X$. This allows us to disprove a conjecture of Marovt (Proc. Amer. Math. Soc. {\bf 134} (2006), 1065-1075). Some related results on other semigroups of functions are also given.
\end{abstract}

\maketitle

\section{Introduction}

For   a compact and Hausdorff  space $X$,  let $C(X, I)$ denote the semigroup of all
continuous functions on $X$ taking values in the closed unit interval $I = [0,1]$. Motivated by a result of L. Moln\'ar (see \cite{Mo}), J. Marovt studied the form of multiplicative bijections on $C(X, I)$
in a special case, namely, when $X$ satisfies the first axiom of countability (see \cite{M}). He proved that
if $\varphi$ is such a map, then there exist a  homeomorphism $\mu : X \ra X$ and a continuous map $p: X \ra (0, +\infty)$ such that
$$\varphi (f) (x) = f(\mu (x))^{p(x)}$$
for every $f \in \dx$ and $x \in X$.
We call maps of this form {\em standard}.  Marovt conjectured  that every multiplicative bijection 
on $\dx$ was necessarily standard, and that the assumption on $X$ of being first countable could be dropped.

In this paper we prove that Marovt's conjecture does not hold. In fact we give a characterization of spaces $X$ for which there exists a multiplicative bijection of $C(X,I)$ that is {\em not}  standard. Roughly speaking, we prove that this is  the case if and only if $X$ is the Stone-\v{C}ech compactification of a proper subset (see Theorem~\ref{sete}).

We study a more general case, as is that of multiplicative bijections between semigroups of functions. Of course, a multiplicative bijection between arbitrary semigroups of $C(X, I)$ need not be of the above form, even if they separate points. A simple example
is the following: take $X$ having only one point, so each semigroup of $C(X, I)$ can be identified with a
semigroup of $I$. Consider for instance $s=1/2$, $t= 1/3$, and $\mathcal{A} := \{s^n t^m : n, m \in \mathbb{N}\}$. It is now clear that the map sending each $s^n t^m$ into $s^m t^n$ is multiplicative and bijective, but cannot be described as above. 

Nevertheless, our technique can be used for semigroups other than $C(X, I)$. In fact, it works if we just  require  these
semigroups $\mathcal{A} \subset C(X, I)$ to satisfy the following three properties (where, as usual,
if $f \in C(X, I)$,  $\mathrm{coz} \hspace{.02in}  f$ denotes the set $\{x \in X: f(x) \neq 0 \}$, and 
 $\mathrm{supp} \hspace{.02in}  f$ denotes its closure).

\begin{description}
\item[Property 1] Given any $x \in X$ and any neighborhood $U$ of $x$,
there exists $f \in \mathcal{A}$ such that $f(x) \equiv 1$ on a neighborhood of $x$ and $\mathrm{supp} \hspace{.02in}  f \subset U$.
\item[Property 2] If $f, g \in \mathcal{A}$ satisfy $f (x) \le  g (x)$ for every $x \in X$,  
and $\mathrm{supp} \hspace{.02in}  f \subset \mathrm{coz} \hspace{.02in}  g $, 
then there exists $k \in \mathcal{A}$ such that
$f=gk$.
\item[Property 3] For each $x \in X$, the set $(0, 1) \cap \{f(x) : f \in \mathcal{A}\} $ is nonempty.
\end{description}

It is easy to see that the image by  a standard multiplicative bijection of a  semigroup satisfying Properties 1, 2, and 3 is a semigroup that also satisfies them. One of such semigroups is, for instance, that of all continuously differentiable functions on $I$ (also taking values in $I$). The general form of multiplicative bijections between this kind of semigroups will be given in Theorem~\ref{derr-laredo}. Finally we also give an example where Marovt's result does not hold for multiplicative bijections defined between general semigroups of this type, even when the ground spaces are first countable (see
Example~\ref{carras-si-cero}). 

\smallskip

Suppose that $\txy$ is multiplicative and bijective, where
 $\cx \subset C(X, I)$ and $\cy \subset C(Y, I)$ are 
semigroups satisfying Properties 1, 2, and 3.  We say that $y \in Y$ is a {\em standard point} for $\varphi$ 
   if there exist a number $p \in (0, + \infty)$ and a point $x \in X$ such that
 $\varphi (f) (y) = f(x)^{p}$ for every $f \in \cx$. 
We denote by $\mathcal{R} (\varphi)$ the set of standard points for $\varphi$.
 
 Throughout we assume that $X$ and $Y$ are compact and Hausdorff spaces. Given a (completely regular) space $Z$,
 we denote by $\beta Z$ its Stone-\v{C}ech compactification.
 
 \smallskip
 
 Remark finally that our results are essentially different from those given in \cite{Mi} for multiplicative bijections between  spaces 
$C(X, \mathbb{R})$ of real-valued continuous functions on $X$. Even so, there is a similarity in that the multiplicative structure of the spaces of functions determines the space
 $X$ up to homeomorphism. More recent results in this line, concerning multiplicative maps on $C(X, \mathbb{R})$, are given for instance in \cite{LS1}  and \cite{LS2}.

 \section{Main results}
 
 We  first give a theorem that provides  a description of multiplicative bijections between semigroups 
 satisfying the above properties. Recall that a subset $Z$ of $Y$ is a {\em cozero-set} if  there exists a real-valued
 continuous function $f$ on $Y$ such that $Z = \{y \in Y: f(y) \neq 0\}$.
 
 \begin{thm}\label{derr-laredo}
 Suppose that $\txy$ is multiplicative and bijective, where
 $\cx \subset C(X, I)$ and $\cy \subset C(Y, I)$ are 
semigroups satisfying Properties 1, 2, and 3. Then 
\begin{enumerate}
\item\label{astamarte} there exist a 
 homeomorphism $\mu$ from $Y$ onto $X$
 and a continuous map $p: \mathcal{R} (\varphi) \ra (0, + \infty)$ such that
$$\varphi(f) (y) = f(\mu(y))^{p(y)}$$
for every $f \in \cx$ and $y \in \mathcal{R} (\varphi)$; 
\item\label{santamarta} the set $\mathcal{R} (\varphi)$ is a  
dense cozero-set in $Y$, and the map $p$ can be extended  to a continuous map from $Y$ into $[0, + \infty]$ taking
values in $\{0 , + \infty\}$ at every point of $Y \setminus \mathcal{R} (\varphi)$. 
\end{enumerate}
\end{thm}

\begin{rem}\label{negro}
Property 3 is necessary to ensure that Theorem~\ref{derr-laredo} holds. For instance if we take $X=Y=[-1,1]$ and 
$\mathcal{A} := \{f \in C(X,I) : f(0) \in \{0,1\}\}$, it is clear that $\mathcal{A}$ satisfies both Properties 1 and 2, and that the map $\varphi: \mathcal{A} \ra \mathcal{A}$ defined, for every $f \in \mathcal{A}$, as $\varphi (f) (x) := f(x)$ if $x \le 0$, and $\varphi (f) (x) := f(x)^2$ if $x >0$, is a multiplicative bijection. Nevertheless, $\mathcal{R} (\varphi) =X$, so each point is standard for $\varphi$, but we cannot find a continuous map $p$ as given in Theorem~\ref{derr-laredo}.
\end{rem}

 We next state the theorem characterizing all spaces on which  can be defined a multiplicative bijection that is not standard. Recall that a topological space $Z$ is said to be {\em pseudocompact} if every real-valued continuous map on $Z$ is bounded.
 
 \begin{thm}\label{sete}
 There exists a bijective and multiplicative  map  $\fxy$ 
 that is not standard if and only if $X$ and $Y$ are homeomorphic and there exists a 
 subset $Z$ of $Y$, not pseudocompact, such that $Y = \beta Z$.
 \end{thm}
 
 Obviously Theorem~\ref{sete} gives an answer in the negative to Marovt's conjecture. It is enough now to take any completely regular space $Z$ (thus ensuring that its Stone-\v{C}ech compactification exists) that is not pseudocompact (as for instance $\mathbb{N}$, $\mathbb{R}$, or any unbounded subset of a normed space), and we have that there are always multiplicative bijections on $C (\beta Z, I)$ that
 are not standard.

\section{Some other results and proofs}

 The following is a key lemma to prove Theorem~\ref{derr-laredo}.

\begin{lem}\label{azero}
Suppose that
$u: \mathscr{A}_1 \ra \mathscr{A}_2$ is an order preserving  multiplicative bijection, 
where $\mathscr{A}_1$ and $ \mathscr{A}_2$ are semigroups contained in $(0,1)$. Then there exists 
$p \in (0 , + \infty)$ such that
$u(\gamma) = \gamma^p$ for every $\gamma \in \mathscr{A}_1$.
\end{lem}

\begin{proof}
Suppose on the contrary that there exist $\alpha, \beta \in \mathscr{A}_1$ such that $u( \alpha) = \alpha^p$ and $u(\beta ) =\beta^q$, where $0<p<q$. By taking integer powers of $\alpha$ if necessary, we may assume without loss of generality that $\alpha < \beta$. Define $a := - \log \alpha$ and $b:= - \log \beta$. Obviously $0 < b <a$, and we can find natural numbers $n, m$ such that 
$$ \frac{b}{a} < \frac{n}{m} <  \frac{q}{p} \frac{b}{a}.$$
Now it is clear that $m b < na$ and $npa <mqb$. These inequalities lead easily to $\alpha^n < \beta^m$ and $\alpha^{np} > \beta^{mq}$, that is, $u(\alpha^n) > u (\beta^m)$, against the fact that $u$ is order preserving.
\end{proof}

\begin{proof}[Proof of Theorem~\ref{derr-laredo}]
(\ref{astamarte})
For each $x  \in X$,  let $\mathscr{U}_x$ be the set of all $f \in \cx$ for which there exists a neighborhood of $x$ where $f \equiv 1 $. 

Let us see that $\mathscr{C}_x := \{y \in Y: \varphi (f) (y) =1 \hspace{.03in} \forall f \in \mathscr{U}_x\}$ is nonempty.
First notice that, given any $f \in \mathscr{U}_x$, the set $\varphi (f)^{-1} (\{1\})$ is compact.
On the other hand, given $f_1, \ldots, f_n \in \mathscr{C}_x$, we can find $f_0 \in \cx$, $f_0 \neq 0$, such that
$\mathrm{supp} \hspace{.02in}  f_0 \subset \{z \in X: \prod_{i=1}^n f_i (z) =1\}$. This obviously implies that
$f_0 f_i =f_0$ for each $i$, so $\varphi(f_0) \varphi(f_i) = \varphi (f_0) \neq 0$. We deduce that 
$\bigcap_{i=1}^n \varphi (f_i)^{-1} (\{1\}) \neq \emptyset$. Then $\{ \varphi (f)^{-1} (\{1\}) : f \in \mathscr{U}_x \}$ satisfies
the finite intersection property, and we conclude that
$\mathscr{C}_x $ is nonempty.

On the other hand,  $\varphi^{-1}$ is also multiplicative, so given any $y \in \mathscr{C}_x$, we  have that the set $\mathscr{C}_y$ (defined in a similar way as $\mathscr{C}_x$) is nonempty. Let us see that $\mathscr{C}_y =\{x\}$. 
Suppose that  there exists  $z \in \mathscr{C}_y$, $z \neq x$. Then we can find $f_z \in \mathscr{U}_x$ such that $z \notin \mathrm{supp} \hspace{.02in}  f_z $.  Since $z \in \mathscr{C}_y$
  and $\varphi (f_z) (y) =1$,
then we can find $k \in \cy$ with $\mathrm{supp} \hspace{.02in} k \subset \mathrm{coz} \hspace{.02in} \varphi (f_z)$
and $\varphi^{-1} (k) (z) \neq 0$. Clearly, 
if we now  take $g \in \mathscr{U}_z$ with $f_z g=0$,  then $g \varphi^{-1} (k) \neq 0$, but $\varphi (g) k =0$, which is impossible.

The above process lets us define a map $\mu : Y \ra X$, which turns out to be bijective, such that $\mu (y)$ is the only point in $\mathscr{C}_y$, and $y$ is the only point
in $\mathscr{C}_{\mu (y)}$, for each $y \in Y$.

We prove next that $\mu$ is continuous at every point of $Y$ (and is consequently a homeomorphism). Take any $y \in Y$, and let $U$ be an open 
neighborhood of $\mu (y)$. We will see that, if $f \in \mathscr{U}_{\mu (y)}$ and $\mathrm{supp} \hspace{.02in}   f \subset U $, 
then $\mu (\mathrm{coz} \hspace{.02in} \varphi (f))$ is contained in $U$. Otherwise
there exists $z \in Y$ such that $\varphi (f) (z) \neq 0$ and $\mu (z) \notin \mathrm{supp} \hspace{.02in}   f$,
so we can take $k \in \mathscr{U}_{\mu (z)}$ such that $kf =0$. Obviously 
$\varphi (k) (z) \varphi (f) (z) \neq 0$, which 
is absurd.

Finally, we have that, by definition, if $y \in \mathcal{R} (\varphi)$, then there exist $p(y) \in (0, + \infty)$
and $x \in X$ such that $\varphi(f) (y) = f (x)^{p(y)}$ for every $f \in \mathcal{A}$. It is easy to check that $x = \mu (y)$.
As for the map $p: \mathcal{R} (\varphi) \ra (0, + \infty)$, we have that for each $y \in \mathcal{R} (\varphi)$, 
 $$p(y) = \frac{\log \varphi (f) (y)}{\log f \pl \mu \pl y \pr \pr}$$for every $f \in \mathcal{A}$ with 
$f \pl \mu \pl y \pr \pr \neq 0,1$. This implies  that $p$ is continuous at $y$, and consequently on $\mathcal{R} (\varphi)$.

\medskip

(\ref{santamarta}) 
For each $y \in Y$,  let $\mathcal{B}_y := \{g (y) : g \in \mathcal{B}\} \cap (0,1)$, and define  $\mathcal{A}_{x}$  for each $x \in X$ in a similar way. Consider also 
 the set $\mathcal{R}_1 (\varphi)$ of all $y \in Y$ such that
$\varphi (f) (y) \neq 0, 1$ whenever $f \in \cx$ satisfies $f (\mu (y)) \neq 0, 1$. We need the following claim.

\smallskip

{\bf Claim.} {\em Let $y \in \mathcal{R}_1 (\varphi)$. If $f, g \in \cx$ satisfy 
$g(\mu (y)) \le f(\mu (y))$,  then $\varphi (g) (y) \le \varphi (f) (y)  $. Moreover $\mu (y) $ belongs to $\mathcal{R}_1 \pl \varphi^{-1} \pr$.}
\smallskip

Suppose first that $g(\mu (y)) < f(\mu (y))$, and  take   a neighborhood $U$ of $\mu (y)$ with $ g(x) < f(x)$ for every $x \in U$. We pick $f_0, g_0 \in \mathscr{U}_{\mu (y)}$
such that $\mathrm{supp} \hspace{.02in}  f_0 \subset U$, and 
such that $\mathrm{supp} \hspace{.02in}  g_0 \subset \{ x \in X : f_0 (x) =1\}$, respectively. Since $\cx$ 
satisfies Property 2, then   
there exists $k \in \cx$ such that $(ff_0)k=  gg_0$. Also $k(\mu(y)) \in (0,1)$, and consequently $$\varphi (g) (y) = \varphi (gg_0) (y) = \varphi (ff_0) (y) \varphi (k) (y) < \varphi (ff_0) (y) = \varphi (f) (y).$$
We now prove that $\mu (y)$ belongs to $\mathcal{R}_1 (\varphi^{-1})$. Let $h \in \cy$ be such that $ h (y) \neq 0,1$.
Suppose that $\varphi^{-1} (h) (\mu (y)) = 0$ and take any $l \in \cx$ with $l (\mu (y)) \neq 0, 1$, and 
$n \in \mathbb{N}$ such that $\pl \varphi(l) (y) \pr^n < h(y)$. We  then have that $ \varphi^{-1} (h)(\mu (y)) < l^n (\mu (y))$ and
$h (y) > \varphi(l^n) (y)$, what goes against what we have proved above. We deduce that $\varphi^{-1} (h) (\mu (y)) \neq 0$. In a similar way we can deduce that $\varphi^{-1} (h) (\mu (y)) \neq 1$. Thus, $\mu (y)$ belongs
to $\mathcal{R}_1 (\varphi)$.

Now, working with $\varphi^{-1}$, it is clear that if $g(\mu (y)) \le f(\mu (y))$, then we cannot get $\varphi (g) (y) > \varphi (f) (y)  $.
The claim is proved.

\smallskip

For any $y \in \mathcal{R}_1 (\varphi)$, we may define a map $\varphi_y : \mathcal{A}_{\mu(y)} \ra \mathcal{B}_y$ in
 the following way. Given $\alpha \in \mathcal{A}_{\mu (y)}$, there exists $f \in \mathcal{A}$ such that
 $f(\mu(y)) = \alpha$. Then define $\varphi_y (\alpha) := \varphi (f) (y)$. It is clear by the above claim that 
 $\varphi_y$ is well defined, and obviously it is multiplicative, order preserving, and bijective. 
Also, we have that  $\varphi (f) (y) =1 $ whenever $f(\mu (y))= 1$, and $\varphi (f) (y) =0 $ whenever $f(\mu (y))= 0$,
$f \in \mathcal{A}$.
Consequently, by 
Lemma~\ref{azero},  we have that $\mathcal{R}_1 (\varphi ) \subset \mathcal{R} (\varphi)$. The other inclusion is immediate, so
$\mathcal{R}_1 (\varphi) = \mathcal{R} (\varphi) $.

\medskip

We  prove that $\mathcal{R} (\varphi) $
is dense in $Y$.
Let $W_0 \subset Y$ be open (and nonempty). Pick 
$y \in W_0$ and assume that  $y \notin \mathcal{R} (\varphi) = \mathcal{R}_1 (\varphi)$, so 
 there exists
$f_0 \in \cx$ such that $f_0 (\mu (y)) \neq 0, 1$ and $ \varphi (f_0) (y) \in \{0,  1\}$. Let $V:= \{x \in X : 0<f_0 (x) <1\}$. 
We next see that $y$ does not belong to the interior of the set $W_1 := \{ z \in W_0 : \varphi (f_0) (z) =0\}$. Otherwise, we can find 
$g_0 \in \mathscr{U}_{y}$ with  $g_0 \varphi (f_0) =0$. This implies that $\varphi^{-1} (g_0) f_0 =0$, against the fact that $\varphi^{-1} (g_0) (\mu (y)) =1$ and $f_0 (\mu (y)) \neq 0$. On the other hand,  
$y$ does not belong to the interior of the set $W_2 := \{ z \in W_0 : \varphi (f_0) (z) =1\}$, because $\varphi (f_0) \notin \mathscr{U}_y$. Consequently, due to the form of $W_1$ and $W_2$, we have that $y$ belongs to the closure of $Y \setminus \pl W_1 \cup W_2 \pr$. This means that 
$$W:= \mu^{-1} \pl V \pr \cap \pl W_0 \setminus \pl W_1 \cup W_2 \pr \pr \neq \emptyset.$$
As in the proof of the claim, it is  clear that given any point in $Y$, and $f , g \in \mathcal{A}$ with $g(\mu (y)) < f(\mu (y))$, then $\varphi (g) (y) \le  \varphi (f) (y)$. Now it is straightforward to check that every point in $W$ belongs to $\mathcal{R}_1 (\varphi) = \mathcal{R} (\varphi)$.

\medskip

We finally see that $\mathcal{R} (\varphi)$ is a cozero-set.  
Suppose  that $y \in Y  \setminus  \mathcal{R} (\varphi)$. Take any net $\pl z_{\alpha} \pr_{\alpha \in \Lambda}$  in $\mathcal{R} (\varphi)$
 converging to $y$.  
Clearly,  since $y \notin \mathcal{R}_1 (\varphi)$, then there exists $f \in \cx$ with $f (\mu (y)) \neq 0, 1$ and 
either $\varphi (f) (y) =0 $ or $1$.
Also  
 $$p(z_{\alpha}) = \frac{\log \varphi (f) (z_{\alpha})}{\log f \pl \mu \pl z_{\alpha} \pr \pr}$$
 for every $\alpha \in \Lambda$. 
  The conclusion that $\lim_{\alpha}  p (z_{\alpha}) \in \{0, + \infty\}$ follows easily.   Finally,  it is clear that
$\mathcal{R} (\varphi)$ coincides with the  cozero-set of the continuous function  $Y \ra \mathbb{R}$ given by 
  $$  y \mapsto p(y) \pl p \pl y \pr ^2 +1 \pr^{-1} $$
 if    $y \in \mathcal{R} (\varphi )$, 
and $y \mapsto 0$ otherwise.  
\end{proof}    

\begin{rem}\label{cerotres}
If in Theorem~\ref{derr-laredo} we assume that the semigroups satisfy just Properties 1 and 2, a similar proof
ensures the existence of the homeomorphism $\mu$ between $Y$ and $X$, even if the theorem does not hold (see Remark~\ref{negro}). 
\end{rem}

In what follows, when we want to specify the homeomorphism $\mu$ and the map $p$ corresponding to a multiplicative and bijective map $\varphi$, as given in Theorem~\ref{derr-laredo}, we will write $\varphi [\mu, p]$.
It is obvious that if $\varphi = \varphi [\mu, p]$, then $\varphi^{-1} = [\mu^{-1}, q]$, 
where $q = 1/\pl p\circ \mu^{-1} \pr$.

\begin{cor}\label{paragogo} 
Let  $\txy$ be multiplicative and bijective,
where
 $\cx \subset C(X, I)$ and $\cy \subset C(Y, I)$ are 
semigroups satisfying Properties 1, 2, and 3. If $\mathcal{R} (\varphi)$ 
is pseudocompact, then $\mathcal{R} (\varphi) = Y$.
\end{cor}

\begin{proof}
It is obvious that if $\varphi = \varphi [\mu, p]$ and $\varphi^{-1} = \varphi^{-1}[ \mu^{-1}, q]$, then both $p$ and $q$ are bounded, and consequently neither $+ \infty$ nor 
$ 0$ are  limit points of $p (\mathcal{R} (\varphi ))$. The conclusion follows from Theorem~\ref{derr-laredo}.
\end{proof}

\begin{prop}\label{dismas}
Let  $\fxy$ be multiplicative and bijective. Then $Y = \beta (\mathcal{R} (\varphi))$.
\end{prop}

\begin{proof}
We follow the notation given  after
 Remark~\ref{cerotres}. 
Suppose that $\varphi = \varphi [\mu, p]$. 
 It is clear that if $\mathbf{1}$ is the map 
constantly equal to $1$ on $X$, and we consider   the multiplicative bijection 
$\psi = \psi[\mu^{-1}, \mathbf{1}] : \dy \ra \dx$,
 then the composition $\varphi \circ \psi : \dy \ra \dy  $ is $\varphi \circ \psi [\mathbf{id}_Y, p] $ (where $\mathbf{id}_Y$  is the identity map on $Y$). Also,  $\mathcal{R}_1 (\varphi) =  \mathcal{R}_1 (\varphi \circ \psi)$, 
so $\mathcal{R} (\varphi) =  \mathcal{R} (\varphi \circ \psi)$ (see Proof of Theorem~\ref{derr-laredo}). Consequently we can assume without loss of 
generality that $X=Y$  and $\mu = \mathbf{id}_Y$. 

Notice  that if we consider the continuous extensions of $p$ and $q$ to maps from $Y$  to $[0, + \infty]$,
 as seen in Theorem~\ref{derr-laredo}, then $p^{-1} (\{0\}) =  q^{-1} (\{ + \infty\}) $
 and $p^{-1} (\{+ \infty\}) =  q^{-1} (\{ 0\}) $. Thus, assuming that $\beta (\mathcal{R} (\varphi)) \neq Y$, at least one of the sets $ p^{-1} (\{0\})$,  $q^{-1} (\{ 0\})$
 is nonempty. We conclude that there is a continuous map 
 from $Y \setminus p^{-1} (\{0\})$, or from $Y \setminus  q^{-1} (\{ 0\})$, to $I$ which does not admit
 a continuous extension to $Y$.
 
 Taking into account that $\mathcal{R} (\varphi) = \mathcal{R} \pl \varphi^{-1} \pr$,  we 
assume without loss of generality that 
  there exists a continuous function $g_0 :Y \setminus p^{-1} (\{0\})  \ra I$ that 
cannot be continuously extended to $Y$. We consider the map $h \in \dy$ whose image by $\varphi$ is the constant function $1/e$. It is easily seen that,
since $p(y) \log h (y) = \log \varphi (h) (y)$ for every $y \in \mathcal{R} (\varphi)$, then 
 $$h (y)= \exp \pl - \frac{1}{p(y)} \pr$$
if $y \in Y \setminus p^{-1} (\{0\})$, and $h \equiv 0$ on $ p^{-1} (\{0\})$.

Since $p^{-1} (\{0\})$ is a zero-set, then there exists a sequence $\pl K_n \pr$ of compact subsets of $Y$ with
$K_{n+1}$ contained in the interior of $K_n$ for each $n \in \mathbb{N}$, and such that 
$p^{-1} (\{0\}) = \bigcap_{n=1}^{\infty}      K_n$. For each $n \in \mathbb{N}$, we take $g_n \in \dy$ satisfying $g_n \equiv g_0$ on
 $Y \setminus K_n$. We denote by $f_n$ its counterimage by $\varphi$.

Notice that if we define $f_0 (z) = f_n (z)$ whenever $z \notin K_n$, then 
$f_0 : 
Y \setminus p^{-1} (\{0\})  \ra I$ is continuous. It is also obvious that $f_0 h$ (defined 
as $0$ on $p^{-1} (\{0\})$) belongs to $C(Y, I)$. On the other hand, we have that $f_0 h \equiv f_n h$ outside $K_n$ for each
$n \in \mathbb{N}$, and consequently $$\varphi (f_0 h) \equiv \varphi (f_n h) \equiv g_n \varphi (h) \equiv g_0 \varphi (h)$$
outside each $K_n$. This implies that $\varphi (f_0 h) \equiv g_0 \varphi (h) $ on $ Y \setminus p^{-1} (\{0\}) $. Now
 it is easy to see that, since $\varphi (h) \equiv 1/e$ on $p^{-1} (\{0\}) $,  then $\varphi (f_0 h)$ is not defined at some points of $ p^{-1} (\{0\}) $, which is absurd.
\end{proof}

\begin{rem}\label{rosales}
We deduce in particular that if the set  $Y \setminus \mathcal{R} (\varphi)$ is nonempty, then its cardinality is at least $2^{\mathfrak{c}}$. Also, if it  is endowed with the restricted topology, then no point in 
$Y \setminus \mathcal{R} (\varphi)$ is a $G_{\delta}$ (see \cite[Chapter 9]{GJ}).
\end{rem}

\begin{rem}\label{azulbelga}
By Remark~\ref{rosales}, we have that each point of $Y$ having a countable base of neighborhoods belongs
to $\mathcal{R} (\varphi)$ for every $\varphi$. In particular, if $Y$ is first countable, then $\mathcal{R} (\varphi) = Y$, which is essentially Marovt's result. But there can be spaces which are not  first countable, and   for which every point is standard. As an easy example,
consider a bijective and multiplicative map $\varphi : C([0, \omega_1], I) \ra C([0, \omega_1], I)$, where 
$\omega_1$ denotes the first noncountable ordinal. Since each point of $[0, \omega_1)$ has a countable base of
neighborhoods, we deduce that $[0, \omega_1] \setminus \mathcal{R} (\varphi) \subset \{\omega_1\}$, and again by Remark~\ref{rosales}, we cannot have $[0, \omega_1] \setminus \mathcal{R} (\varphi) = \{\omega_1\}$. We conclude that $\mathcal{R} (\varphi) = [0, \omega_1]$.
\end{rem} 

\begin{rem}\label{nikdezir}
The conclusion given in Remark~\ref{azulbelga} is not true for more general semigroups, that is, not every point having a countable base of neighborhoods necessarily belongs to $\mathcal{R} (\varphi)$ (see Example~\ref{carras-si-cero}).
\end{rem}

Using Theorem~\ref{derr-laredo}, it is easy to see that every isolated point is standard for every multiplicative bijection. More generally, the next result allows us to identify some points that belong to $\mathcal{R} (\varphi)$ for every $\varphi$.
    
\begin{cor}\label{afeite}
Let $y \in Y$ be such that $\beta \pl Y \setminus \{y\} \pr \neq Y$. 
Then $y \in \mathcal{R} (\varphi)$ for every multiplicative bijection $\fxy$.
\end{cor}

\begin{proof}
Suppose that $y \notin \mathcal{R} (\varphi)$. Thus, taking into account  Proposition~\ref{dismas} and the fact that $\mathcal{R} (\varphi) \subset  Y \setminus \{y\}$, we deduce that $\beta \pl Y \setminus \{y\} \pr = Y$, and we 
are finished.
\end{proof}

\begin{rem}
Obviously, the converse of Corollary~\ref{afeite} is in general not true. For an easy example, consider for instance
 the point $\omega_1$ in Remark~\ref{azulbelga}, and take into account that  $\beta [0, \omega_1) = [0, \omega_1]$.
\end{rem}

\begin{rem}
Marovt's result cannot  be given in general for the kind of semigroups we are dealing with, even if the ground spaces are assumed to be first countable. That is, it is possible that $\mathcal{R} (\varphi) \neq X$ for a multiplicative bijection $\varphi : \cx \ra \cy$, even in the case when $\mathcal{A}, \mathcal{B} \subset \dx$ satisfy Properties 1, 2, and 3, and  the 
space $X$ is first countable. We next include an example of this fact.
\end{rem}

\begin{ex}\label{carras-si-cero}
For a topological space $Z$, denote by $C(Z, \ro)$ the set of all continuous maps from $Z$ to $[0, +\infty)$. Let $\nin $ be the one-point compactification of $\mathbb{N}$ and  let 
$\varphi : C( \nin , \ro ) \ra C ( \mathbb{N} , \ro )$ be defined, for each $f \in C( \nin , \ro) $, as $\varphi (f) (n) := f(n)^n$ for every $n \in \mathbb{N}$.
 It is clear that $\varphi$ is multiplicative and injective. 
Since it also preserves order, then the image of each $f \in C(\nin , I)$ takes values in $I$.

Consider now  the subset $\mathcal{A}_1$ of $C(\nin, I )$
of all characteristic  functions which are continuous, and let $f_0 \in C( \nin , I)$   be the counterimage by $\varphi$ of the constant function equal to $1/2$. Define $\mathcal{A}$ as the set of all functions of the form
$f = \alpha g f_0^m$ for some $\alpha \ge 0$, $g \in \mathcal{A}_1$, and $m \in \mathbb{Z}$, such that $f(n) \in I$ for every $n \in \mathbb{N}$.

 It is easy to see that if $f \in \mathcal {A}$, then $\lim_{n \ra \infty} \varphi (f) (n)$ exists, so in a natural way
 we may define a map $\varphi : \mathcal{A} \ra C( \nin, I) $.
Put  $\mathcal{B} := \varphi (\mathcal{A})$. It is easy to see that 
$\mathcal{A}$ and $\mathcal{B}$ satisfy Properties 1, 2, and 3. On the other hand, every point in $\mathbb{N}$ is
 isolated, and consequently $\mathbb{N}$ is contained in $\mathcal{R} (\varphi)$. Nevertheless, $\infty$ does not 
belong to $\mathcal{R} (\varphi)$.
\end{ex}

\begin{proof}[Proof of Theorem~\ref{sete}]
Suppose  that there exists a multiplicative bijection $\fxy$ that is not standard. First, by 
Theorem~\ref{derr-laredo}, $X$ and $Y$ must be homeomorphic. Also, calling $Z := \mathcal{R} (\varphi)$, we have 
by Corollary~\ref{paragogo} that $Z$ is not pseudocompact, and 
by Proposition~\ref{dismas} that
$Y = \beta Z$.

Conversely, suppose that $X$ and $Y$ are homeomorphic, and that there exists a proper subset $Z$ of $Y$ which is not pseudocompact, and such that
$Y= \beta Z$. It is clear that without loss of generality we may assume that $X=Y$. 
Since $Z$ is not pseudocompact,  then there exists an unbounded continuous function $u: Z \ra [1 , + \infty)$. We define a  map $\varphi: \dz \ra \dz$ as $\varphi (g) (z) := g (z)^{u(z)}$ for each $g \in \dz$ and $z \in Z$. 

It is easy to check that $\varphi$ is multiplicative and injective.
Let us see that $\varphi$ is surjective. To this end, we take any  $k \in \dz$, and consider the map
 $L_k : Z \ra [- \infty, 0]$
defined as $L_k (z) := (\log k(z)) / u(z)$ if $k (z) \neq 0$, and $L_k (z) := -\infty$ if $k(z) =0$. It is easy to 
check that $L_k$ is continuous and that $\varphi (\exp L_k) = k$ (assuming $\exp (- \infty) =0$).

Obviously $\varphi$ can be seen as a multiplicative bijection on $C(\beta Z, I)$. On the other hand,  if 
for some homeomorphism $\mu: \beta Z \ra \beta Z$ and a continuous map 
$v: \beta Z \ra (0, + \infty)$,  
$\varphi (f) (x) = f(\mu (x))^{v (x)}$  for every $f$, then $\mu$ must be the identity on $\beta Z$, and $v(z) = u(z)$ for every $z \in Z$. 
This obviously implies that $v$ must attain the value $+ \infty$ on some points of $\beta Z \setminus Z$, which is 
absurd.
\end{proof}

\section{Acknowledgements}
The author wishes to thank  the referee for his/her  valuable remarks and, also,  Ana M. R\'odenas for drawing his attention to this subject.

\end{document}